\documentclass[12pt]{article}
\usepackage[dvips]{graphicx}
\usepackage{latexsym}
\usepackage{amsfonts}
\newtheorem{coro}{Corollary}[section]
\newtheorem{defi}{Definition}[section]
\newtheorem{lema}{Lemma}[section]
\newtheorem{rema}{Remark}[section]
\newtheorem{notac}{Notation}[section]
\newtheorem{ejem}{Example}[section]

\newtheorem{question}{Question}[section]
\newtheorem{prop}{Proposition}[section]
\newtheorem{teor}{Theorem}[section]

\title{\bf Helix, shadow boundary and minimal submanifolds}
\author{}
\date{}
\begin{document}
\maketitle

\begin{abstract}
Inspired by a Blaschke's work about analytic convex surfaces, we
study {\em shadow boundaries} of Riemannian submanifolds $M$,
which are defined by a parallel vector field along $M$. Since a
shadow boundary is just a closed subset of $M$, first, we will
give a condition that guarantee its smoothness.
It depends on the second fundamental form of the submanifold.
It is natural to search for what kind of properties might have such
submanifolds of $M$? Could they be totally geodesic or minimal?
Answers to these and related questions are given in this work.
\end{abstract}

\noindent {\footnotesize {\bf Key words:} minimal and helix
submanifolds, smooth shadow boundaries, parallel vector field,
holonomy group.}

\begin{section}{Introduction}

In their book \cite{No-Sa}, Nomizu and Sasaki
defined the shadow boundary of a Euclidean surface $S$ with respect
to a fixed direction $v$ as the set of points $p \in S$ such that: the line through $p$ in the direction of $v$ is tangent to $S$
(see also \cite{Bu}).
Blaschke characterized convex analytic surfaces
with planar shadow  boundaries. We can see an extension
of his result in \cite{No-Sa}, page 61.\\
Recently M. Ghomi in \cite{Gho}, solved the shadow problem
formulated by H. Wente. The shadow of a Euclidean orientable
surface $S$, with unit normal vector field
$n: S \longrightarrow S^2$, consist of those
points $p \in S$ where the inner product between $n(p)$ and $v$ is
positive. His result says that a closed surface with simply
connected shadows should be convex (see also \cite{Ho2}).
He used {\em horizon} term in instead of shadow boundary.\\

In a previous work \cite{GRH2}, we studied shadow
boundaries of any Euclidean submanifold.
Now, we extend the concept of shadow boundary in the
following way. Let $N$ be a Riemannian manifold and let
$M \subset N$ be a submanifold. Let us assume that
$Y:M \longrightarrow TN$ is an invariant vector field under
parallel transport in $N$ along curves contained in $M$,
in Theorem \ref{existenciaCampoParalelo-cr-subvariedad}
we give a condition in terms of the holonomy group of $N$
for the existence of such vector field $Y$.
For example $Y$ could be the restriction to $M$ of a global
parallel vector field on $N$.
The {\em shadow boundary} of $M$ {\em with respect to} $Y$
is the following subset of $M$.
$$S\partial(M,Y)=\{ x \in M\ | \ Y(x) \in T_xM \},$$
i.e. the points where the vector field is tangent to the
submanifold. This is a closed subset of $M$. In Theorem
\ref{Frontera-de-sombra-suave}, there is a
condition over the second fundamental form of $M \subset N$
that makes $S\partial(M,Y)$ a submanifold of $M$. A
consequence is that $S\partial(M,Y)$ could be a finite set
of points, if $M$ is compact with $\dim N= 2 \dim M$.\\

When we were looking for properties of these type of submanifolds
of $M$, we observed the necessity to introduce Helix
submanifolds of a Riemannian manifold: Submanifolds which make
constant angle with $Y$ (see definition \ref{helix}). We proved in
Theorem \ref{tgs-implica-helice} that if $L \subset
S\partial(M,Y)$ is a totally geodesic submanifold of $M$, then it
is a helix submanifold of $N$ with respect to $Y$. A special case
of helix submanifold is when the constant angle is orthogonal. In
theorem \ref{ortogonal-implica-tgs} we proved that if $L \subset
S\partial(M,Y)$ is orthogonal helix of $N$,
then $L$ is a totally geodesic submanifold of $M$.\\

Theorem \ref{helix-codimension-one} might be of help to develop
some intuition about helix submanifolds of codimension one.
Theorem \ref{coro-helice-orto-stg} implies that a submanifold of $N$
which is helix with respect to cod$-M$ parallel vector  fields is
a totally geodesic submanifold of $N$. Also Corollary
\ref{torus-helix} proves that a closed orientable helix surface
in a three-dimensional manifold should
be a torus or a totally geodesic submanifold.\\

Minimal submanifolds are relatives of totally geodesic
submanifolds. A submanifold contained in a shadow boundary,
$L \subset S\partial(M,Y)$, might be a minimal submanifold of
$M$: Theorem \ref{minimal-contenida-fs} tell us that a
necessary and sufficient condition is that the mean curvature
vector field of $L \subset N$ should be orthogonal to $Y$.
Finally, if $N=\mathbb{R}^{n+1}$ and $L \subset S\partial(M,Y)$
is compact, minimal of codimension one in $M$ and $L$ is contained
in a totally geodesic submanifold of $N$, then $L$ is totally geodesic in $M$.

\begin{notac}
\em
We will work on the $C^\infty$ category.
In this manuscript $(N,g)$ will be a Riemannian
manifold with covariant derivative $\nabla$.
We will use $\mathfrak{X}(N)$ to denote vector
fields on $N$.
\end{notac}

\end{section}

\begin{section}{Parallel vector fields along submanifolds}

Let $M \subset N$ be a submanifold.
For every $x \in M$, there is a decomposition:
$$T_xN=T_xM \oplus T_xM^\perp.$$
This is a direct sum, i.e., there is a unique
decomposition for every $V \in T_xN$:
$V=\mbox{tan}(V) + \mbox{nor}(V)$. Where
$\mbox{tan}(V) \in T_xM$ and $\mbox{nor}(V) \in T_xM^\perp$.
From this, we can define two natural applications,
$\mbox{tan}:TN \longrightarrow TM$ and
$\mbox{nor}: TN \longrightarrow TM^\perp$.\\
So, every vector field $Y: M \longrightarrow TN$ along $M$
may be decomposed as $Y=\mbox{tan}(Y)+\mbox{nor}(Y)$ into two
vector fields (see \cite{O} for details).

\begin{notac}
\em
We will use $\mathfrak{X}(N)$ to denote the vector fields on
$N$. Given a submanifold $M$ of $N$, we can consider
vector fields in $N$ along $M$:\\
$\mathfrak{X}(N,M)=\{ Y: M \longrightarrow TN \}$.
Each $Y \in \mathfrak{X}(N,M)$ induces two natural vector fields
$\mbox{tan}(Y): M \longrightarrow TM$ and
$\mbox{nor}(Y):M \longrightarrow TM^\perp$.

\end{notac}

Let us remember the classic definition of a parallel
vector field on a manifold.

\begin{defi}
\em
A vector field $X \in \mathfrak{X}(N)$
is {\em parallel} if $\nabla_W X=0$ for every
$W \in TN$.
\end{defi}

\begin{ejem}
\em
Let $N= \mathbb{R} \times M$ be a Riemannian product manifold.
Then $N$ admits a parallel vector field defined by
$Y: N \longrightarrow TN$ with
$Y(t,y)=\partial_t$, where
$T_{(t,y)}N = T_t\mathbb{R} \oplus T_yM$.
\end{ejem}

\begin{rema}
\em
In \cite{WeI} and \cite{WeII}, D. J. Welsh  studied
the existence of a parallel vector field on a Riemannian
manifold. In \cite{WeI} he gives the following criterion.\\
A complete and connected Riemannian manifold $N$ admits
$p$ linearly independent parallel vector fields
if and only if there exists a Riemannian manifold $M_2$,
and a group $L \subset \mathbb{R}^p \times
I(M_2)$ such that\\
(a) the first projection $pr|L$ is injective and\\
(b) the orbits of $L$ in $\mathbb{R}^p \times M_2$ are discrete,
so $N$ is isometric to $(\mathbb{R}^p \times M_2)/L$.
\end{rema}

\begin{defi}
\label{parallel-along-submanifold}
\em
Let $M$ be a Riemannian submanifold of $N$ and let
$Y \in \mathfrak{X}(N, M)$. We will say that $Y$ is a
{\em parallel vector field along} $M$, if $\nabla_W Y=0$
for every $W \in TM$.
\end{defi}

\begin{notac}
\em
We will denote by $\mathfrak{X}_0(N,M)$ the set of all vector
fields $Y: M \longrightarrow TN $ which are parallel
along $M$.
\end{notac}

\begin{ejem}
\em
Let us assume that $N$ admits a (global) parallel vector field
$Y \in \mathfrak{X}(N)$. If $M$ is a submanifold of $M$ then
$Y_{|M}: M \longrightarrow TN$ is a parallel vector field
along $M$.\\
Let us consider another example. Let $N$ be a submanifold
of $\mathbb{R}^n$ and let $H$ be any linear submanifold
of $\mathbb{R}^n$ tangent to $N$. Assume that $M=N \cap H$
is a submanifold. Then any constant vector field $X$ in
$H$ induces a parallel vector field along
$M$, $Y=X_{|M}: M \longrightarrow TN$, where $Y$ it is not
necessarily the restriction of a parallel vector field on $N$.
\end{ejem}

\begin{rema}
\em Let $M \subset N$ be a submanifold of codimension one. Let
assume that $Y: M \longrightarrow TN$ is a parallel vector field
along $M$. If $Y$ tangent to $M$ (i.e. $Y \in \mathfrak{X}(M)$),
then $Y$ is a parallel vector field on $M$. If $Y$ is orthogonal
to $M$, the it is parallel with respect to the normal connection
$\nabla^\perp$ on $TM^\perp$. Let us observe that the converse
assertions are false. But if $M$ is a totally geodesic
submanifold, then a parallel vector field on $M$ (connection on
$M$) or normal to $M$ (normal connection) is parallel along $M$
(connection on $N$).
\end{rema}

\begin{rema}
\em Let $M \subset N$ be a Riemannian submanifold with induced
connection $\nabla^M$. Let $II_x: T_xM \times T_xM \longrightarrow
T_xM^\perp$ be the second fundamental form of $M \subset N$ at $x
\in M$. Finally, let $\nabla^\perp$ the
induced normal connection on $TM^\perp$.\\
If $X :M \longrightarrow TM$ is a parallel vector field on $M$
($\nabla^M X=0$) and for every $x \in M$, $II_x(X(x), \cdot)=0$,
then $X$ is parallel along $M$, i.e. $\nabla_W X=0$
($W \in T_xM$).\\
If $Z :M \longrightarrow TM^\perp$ satisfies $\nabla^\perp Z=0$
and for every $x \in M$, $g(Z(x),II_x(\cdot, \cdot))=0$, then
$Z$ is parallel along $M$, i.e. $\nabla_W Z=0$ ($W \in T_xM$).
\end{rema}

\begin{rema}
\em
Let $M$ be a connected submanifold of $N$.
If $Y$ is a parallel vector field along $M$, then
$g(Y,Y)$ is constant:
Let $\alpha$ be a smooth curve in $M$ trough $p \in M$.\\
So, $\frac{d}{dt}g(Y(\alpha(t)), Y(\alpha(t)))= 2g(Y(\alpha),
\nabla_{\dot\alpha} Y(\alpha))=0$ because $Y$ is parallel along
$M$. Therefore $g(Y,Y)=g(Y(p),Y(p))$ and since $M$ is connected,
there exist a smooth curve in $M$ from any point to $p$.
\end{rema}

\begin{defi}
\label{holonomia-relativa} \em Let $M$ be a Riemannian submanifold
of $N$. For every $y \in N$, we denote by $Hol_y(N)$ the Holonomy
group based at $y$ of the Levi-Civita connection of $(N,g)$. This
is a subgroup of $O(T_yN)$, i.e. its elements are isometries of
$T_yN$. To describe $Hol_y(N)$, we should consider a loop $\gamma$
based at $y$ (piece-wise smooth closed path through $y$) and take
the parallel transport $P_\gamma$ along it. Then
$$Hol_y(N)=\{P_\gamma : T_yN \longrightarrow T_yN \ | \
\gamma \subset N \mbox{ loop based at } y \}.$$
Let us consider the following subgroup of $Hol_x(N)$,
where $x \in M$,
$$Hol_x(N,M)=\{P_\gamma \in Hol_x(N) \ | \ \gamma \subset M  \}.$$
We will call $Hol_x(N,M)$, the {\em Holonomy
subgroup at $x$ with respect} to $M$.\\
Let us remember that $Hol_y(N)$ acts in $T_yN$.
\end{defi}
\begin{teor}
\label{existenciaCampoParalelo-cr-subvariedad} Let $M \subset N$
be a submanifold and let $G=Hol_x(N,M)$, where $x\in M$. There
exist a parallel vector field along $M$, $X:M \longrightarrow TN$,
if and only if there exist $W \in T_xN$ such that $G_W=G$, where
$G_W$ is the isotropy subgroup at $x$ under the action of
$G$ in $T_xN$.
\end{teor}
{\bf Proof.}
Let $W \in T_xN$ be a fixed vector
under the action of $G$ in $T_xN$.
For every $z \in M$, we define $X(z)=P_\beta(W)$, where
$\beta : [0,1] \longrightarrow M$ is any piece-wise smooth
regular curve with $x=\beta(0)$ and $z=\beta(1)$.
Affirmation: $X$ does not depends on $\beta$.
Let $\alpha$ be another curve with the same conditions as
$\beta$.
Let us consider the next loop
$$
\gamma(t)=
\left\{
\begin{array}{ccc}
\tilde{\beta}(t)=\beta(2t) & \mbox{if} & t \in [0,\frac{1}{2}]\\
\tilde{\alpha}(t)=\alpha(2-2t) & \mbox{if} & t \in [\frac{1}{2},1]
\end{array}
\right.
$$
By hypothesis,
$P_\gamma(W)=W$, i.e., $P_{\tilde{\alpha}}
(P_{\tilde{\beta}}(W))=W$. It follows that
$P_{\tilde{\beta}}(W)=P^{-1}_{\tilde{\alpha}}(W)$.
Then $P_\beta(W)=P_\alpha(W)$.
Here we used that we have a homomorphism group from
$\pi_1(N,x)$ into $Hol_x(N)/Hol_0(N)$, where $Hol_0(N)$ is
the restricted Holonomy group
(null-homotopic loops, see \cite{Be}, page 280).\\
The smoothness of $X$ follows from observations in Besse's
book \cite{Be}, page 282.
$\Box$

\begin{defi}
\em
Let $M \subset N$ be a  Riemannian submanifold. If
every geodesic of $M$ is a geodesic of $N$, then
$M$ is called a {\em totally geodesic submanifold}
({\bf tgs}).
\end{defi}

\begin{teor}
\label{coro-helice-orto-stg}
Let $M \subset N$ be a Riemannian submanifold
of codimension $r$. Let $X_j: M \longrightarrow TN$,
$j=1,\ldots,r$, be parallel vector fields along $M$,
such that for every $x \in M, \ \ \{X_1(x),
\ldots, X_r(x)\}$ is a basis of $T_xM^\perp$.
Then $M$ is a totally geodesic submanifold of $N$.
\end{teor}
\noindent {\bf Proof.}
Let $x \in M$, by hypothesis,
$\dim(T_xM^{\perp}) =r$, and
\begin{equation}
\label{descom-mayor-dim}
T_xN= T_xM \oplus (\oplus_{j=1}^r <X_j(x)>).
\end{equation}
Let $\gamma \subset M$ be a geodesic through $x$,
so $\nabla_{\dot{\gamma}} \dot{\gamma}$
is orthogonal to $M$, i.e.,
$\nabla_{\dot{\gamma}} \dot{\gamma} \in {T_\gamma M}^\perp.$
Affirmation: $\gamma$ is geodesic of $N$.
By equality (\ref{descom-mayor-dim}),
we need prove that $\nabla_{\dot{\gamma}} \dot{\gamma}$
is orthogonal to every $X_j(\gamma)$:
Since $g(X_j(\gamma), \dot{\gamma})=0$,
$$g(X_j(\gamma), \nabla_{\dot{\gamma}}\dot{\gamma})+
g(\nabla_{\dot{\gamma}}X(\gamma), \dot{\gamma}) =
dg(X_j(\gamma), \dot{\gamma})=0.$$
Therefore,
$\nabla_{\dot{\gamma}} \dot{\gamma}$ is orthogonal to $N$.
Hence, $\gamma$ is a geodesic of $N$.
$\Box$

\end{section}

\begin{section}{Helix submanifolds}

The next definition is a natural extension of the classic
concept of general helix in $\mathbb{R}^3$ which appears from
the first courses in differential geometry: a curve in
$\mathbb{R}^3$ which makes constant angle with respect to a fixed
direction. There are extensions into a  three manifold, but the
helix is again a curve (see \cite{Ba} and \cite{E-I} ). In the following definition a helix may be
a submanifold of higher dimension.

\begin{defi}
\label{helix}
\em
Let $M$ be a Riemannian submanifold of $N$
and let $Y \in \mathfrak{X}_0(N, M)$ be a parallel vector
field along $M$. We say that $M$ is a {\em helix submanifold}
of $N$ {\em with respect} to $Y$ if the following function
$h: M \longrightarrow \mathbb{R}$ is constant.
\begin{equation}
\label{funcion-angulo}
h(x)=\mbox{max}\{g(w, Y(x)) | w \in T_xM, \ g(w,w)=1 \}.
\end{equation}
\end{defi}

\begin{rema}
\label{angulomaximo}
\em
We could think of $h(x)$ as the angle between $T_xM$ and
$Y(x)$. Let us observe that $h(x)=g( \frac{\mbox{tan}(Y(x))}
{(g(\mbox{tan}(Y(x)), \mbox{tan}(Y(x))))^{1/2}}, Y(x) )
=(g(\mbox{tan}(Y(x)), \mbox{tan}(Y(x))))^{1/2}$.
\end{rema}

\begin{ejem}
\em
\label{totallygeodesichelix}
Let $M \subset N$ be a connected and totally geodesic submanifold.
If $Y \in \mathfrak{X}_0(N,M)$ then $M$ is a helix submanifold
of $N$ with respect to $Y$. By Theorem
\ref{existenciaCampoParalelo-cr-subvariedad},
$Y$ is invariant under parallel transport on $N$
along curves contained on $M$. Is well known that $TM$ is also invariant under parallel transport on $N$.
So the angle between $Y$ and $TM$ is constant.
\end{ejem}

\begin{ejem}
\em
Let $\gamma \subset N$ be a embedded geodesic.
Since the tangent vector $\dot\gamma$ of $\gamma$ is parallel
along it, $\gamma$ is a helix of $N$ with respect to any
parallel vector field ($\dot\gamma$ itself) along $\gamma$.\\
A helix submanifold is not necessarily a totally geodesic
submanifold, for example a circular cylinder and any cone of  revolution in $\mathbb{R}^3$ are helix submanifolds with
respect to a constant  vector field parallel to their axis.
\end{ejem}

\begin{lema}
\label{proy-constante}
Let $M$ be a connected helix submanifold of $N$ with respect
to $Y \in \mathfrak{X}_0(N,M)$. Then $\mbox{\em tan}(Y)$ and
$\mbox{\em nor}(Y)$ have constant lenght.
\end{lema}
{\bf Proof.}
Since $M$ is connected, $Y$ has constant lenght.\\
By hypothesis, the function
$h$ in (\ref{funcion-angulo}) is constant, so remark
\ref{angulomaximo} implies that
$h(x)=(g(Y_0(x), Y_0(x)))^{1/2}$ is constant, where
$Y_0(x)=\mbox{tan}(Y(x))$.
To prove that $\mbox{nor}(Y)$ has constant lenght, we
can use the decomposition $Y=\mbox{tan}(Y)+\mbox{nor}(Y)$
and the equality $g(Y,Y)=g(\mbox{tan}(Y), \mbox{tan}(Y))+
g(\mbox{nor}(Y), \mbox{nor}(Y))$.
$\Box$

\begin{coro}
\label{zeroEuler}
Let $M$ be a compact, orientable and connected submanifold of
$N$. If $M$ is a helix submanifold with respect to $Y$, then
it has zero Euler characteristic or $Y$ is orthogonal to $M$.
\end{coro}
{\bf Proof.}
Let us assume that $Y$ is not orthogonal to $M$. Then by Lemma
\ref{proy-constante}, $\mbox{tan}(Y) \in \mathfrak{X}(M)$ is
a vector field on $M$ of constant lenght. Since $M$ is compact and  orientable we conclude, by a well known Poincare-Hopf's theorem
(see \cite{GuiPo}), that $M$ has zero Euler characteristic.
$\Box$

\begin{coro}
\label{torus-helix}
Let $M$ be a connected, orientable and compact submanifold
of dimension two of $N^3$. If $M$ is a helix of $N$, then $M$
is diffeomorphic to a torus or it is totally geodesic.
\end{coro}
{\bf Proof.}
Let assume that $M$ is a helix submanifold with respect to
$Y \in \mathfrak{X}_0(N,M)$, a parallel vector field along $M$.
By Corollary \ref{zeroEuler}, $M$ has zero Euler characteristic or $Y$
is orthogonal to $Y$. In the first case, we can deduce that
$M$ is a torus. In the second case, $M$ has an orthogonal
vector field. Since $Y$ is parallel along $M$, then $M$ is totally
geodesic.
$\Box$
\begin{ejem}
\label{Euclidean-helix}
\em
Let us consider a connected hypersurface $M$ in
$N=\mathbb{R}^{n+1}$.
Let $Y \in \mathfrak{X}(N)$ be a constant vector field.
If $M$ is a helix submanifold of $N$, then
\begin{itemize}
\item $M$ is contained in a hyperplane (perpendicular to $Y$) of
      $N$ when $Y$ is orthogonal to $M$.
\item $M$ is not compact (otherwise $Y$ would be orthogonal to
      $M$),
\item $M$ is orientable ($\mbox{nor}(Y)$ induces an orientation),
\item $M$ has zero Gauss-Kronecker curvature (the Gauss map of
      $M$ is singular).
\end{itemize}
If $M \subset \mathbb{R}^{n+1}$ is not a hypersurface but
is compact, we can conclude that $M$ is contained in a
hyperplane orthogonal to $Y$.
\end{ejem}

\begin{prop}
If $M$ is a compact helix of $N=\mathbb{R} \times M_2$ with
respect to $X=\partial_t$ then $X$ is orthogonal to $M$.
\end{prop}
{\bf Proof.}
Since $M$ is compact, the projection $\pi_1$ of
$M$ in $\mathbb{R}$ is compact, so $\pi_1(M)
\subset \mathbb{R}$ has a maximum denoted by
$t_0$. Let $x \in M$ be such that $\pi_1(x)=t_0$,
affirmation: ${t_0} \times M_2 =\pi_1^{-1}(t_0)$
is tangent to $M$ in $x$. We deduce from this that
$T_xM \subset T_x ({t_0} \times M_2)$.
Let us observe that $X$ is orthogonal to ${t_0} \times M_2$,
in consequence $X$ is orthogonal to $M$
at $x$. Since $M$ is a helix, $X$ is orthogonal to $M$.
$\Box$

\begin{rema}
\em
In general a compact helix submanifold $M$, with respect to a
global parallel vector field $X$ on $N$, is not necessarily orthogonal to $X$. Using Example \ref{totallygeodesichelix},
we can construct an example of this: $M$ should be a totally
geodesic submanifold of $N$ and $X$ a global parallel vector field
on $N$ but non-orthogonal, for example tangent along $M$.
\end{rema}

For general Riemannian hypersurfaces which are helix, we can
prove the following result.
\begin{teor}
\label{helix-codimension-one}
Let $M$ be a connected hypersurface in a
Riemannian manifold $N$.
Let us assume that $M$ is a helix submanifold of $N$ with
respect to $Y \in \mathfrak{X}_0(N,M)$. Then\\
a) If $Y$ is orthogonal to $M$ at some point, then\\ $M$ is
   totally geodesic submanifold of $N$.\\
b) If $Y$ is tangent to $M$ at some point, then\\ $M$ is locally
   a Riemannian product $\mathbb{R} \times M_2$.\\
c) If $Y$ is transversal to $M$ at some point and the integral
   curves of $\mbox{tan}(Y)$ are geodesics in $M$, then they
   are geodesics in $N$.
\end{teor}
{\bf Proof.}
Since $M$ is a helix, the angle between $Y$ and $M$ is constant.\\
a). We have that $Y$ is parallel along $M$ and orthogonal to $M$.
So, $M$ is a totally geodesic submanifold of $N$.\\
b). Let us observe that $Y$ is a parallel vector field on
$M$, then by Welsh's work in \cite{WeI}, $M$ is locally
isometric to a Riemannian product.\\
c). In this case, $Y$ is transversal to $M$ in any point.
Let $Y_0=\mbox{tan}(Y)$ and let $\alpha \subset M$ be an integral
curve of $Y_0$, i.e. $\dot\alpha (t)= Y_0(\alpha (t))$. Affirmation:
$Y$ is orthogonal to $\nabla_{\dot{\alpha}}\dot{\alpha}$.
$$0=\frac{d}{dt}g(\dot{\alpha}, Y_0(\alpha))=
\frac{d}{dt}g(\dot{\alpha}, Y(\alpha))=
g(\nabla_{\dot{\alpha}}\dot{\alpha}, Y)+
g(\dot{\alpha}, \nabla_{\dot{\alpha}} Y)=g(\nabla_{\dot{\alpha}}
\dot{\alpha}, Y),$$
where $\nabla_{\dot{\alpha}} Y =0$ because $Y$ is parallel along
$M$. Since $\alpha$ is geodesic in $M$, $\nabla_{\dot{\alpha}}
\dot{\alpha}$ is orthogonal to $M$. So, $\nabla_{\dot{\alpha}}
\dot{\alpha}=0$, otherwise $Y$ would be tangent to $M$ but it is
transversal to $M$.
$\Box$

\begin{rema}
\em
Let $M \subset N$ be a submanifold of codimension one and let
$Y \in \mathfrak{X}_0(N,M)$ be transverse (important) to $M$.
If every geodesic of $M$ is a helix of $N$ with respect to $Y$,
then $M$ is a totally geodesic submanifold of $N$. Proof:
Let $\gamma$ be a geodesic on $M$
($\nabla_{\dot{\gamma}}\dot{\gamma} \in TM^\perp$), then
the equation
$$
g(Y(\gamma(t)), \nabla_{\dot{\gamma}}\dot{\gamma})+
g(\nabla_{\dot{\gamma}}Y(\gamma(t)), \dot{\gamma}) =
\frac{d}{dt}g(Y(\gamma(t)), \dot{\gamma})$$
and the hypothesis imply that
$g(Y(\gamma(t)), \nabla_{\dot{\gamma}}\dot{\gamma})=0$.
Since $M$ is of codimension one ,
$T_{\gamma(t)}N=<Y(\gamma(t))>\oplus T_{\gamma(t)}M$.
Therefore,  $\nabla_{\dot{\gamma}}\dot{\gamma}=0$,
i.e. $\gamma$ is a geodesic of $N$.\\
If $Y$ is not transverse to $M$, we can not conclude that
$M$ is a totally geodesic submanifold of $N$. For example
a circular cylinder $M$ in $N=\mathbb{R}^3$. Their geodesics are
helix (and $M$ itself) with respect to a constant vector field
in the direction of its axis.
\end{rema}

\begin{rema}
\label{minimalen-dim-tres}
\em
Let $N=\mathbb{R}^3$ and let $v \in N$ be a non-zero vector.
If $M^2 \subset N$ is a complete minimal surface which
is a helix with respect to $Y=v$, then $M$ is a plane.\\
Proof: By Example \ref{Euclidean-helix}, $M$ has zero
Gauss-Kronecker curvature, but in dimension
$2$ it is the Gaussian curvature. Since $M$ is minimal,
it easy to see that $M$ is a plane.
\end{rema}

\begin{question}
\em
If $M^n \subset N^{n+1}$ is minimal and a helix submanifold,
does it is a totally geodesic submanifold of $N$?
\end{question}

\begin{rema}
\em The argument in Remark \ref{minimalen-dim-tres}, when
$N=\mathbb{R}^3$, is not valid in $N=\mathbb{R}^{n+1}$ with $n
\geq 3$: In $\mathbb{R}^4$ there are minimal hypersurfaces with
zero Gauss-Kronecker curvature, like $M^3=M' \times
\mathbb{R}$, where $M'$ is minimal in $\mathbb{R}^3$.
\end{rema}

\end{section}

\begin{section}{Shadow boundary and helix}

\begin{defi}
\em
Let $M$ be a Riemannian submanifold of $N$, and
let $Y : M \longrightarrow TN$ be a parallel
vector field along $M$
(i.e. $Y \in \mathfrak{X}_0(N,M)$).
The {\em shadow boundary} of
$M$ {\em with respect to} $Y$ is the following subset
of $M$.
\begin{equation}
\label{sb-en-campos-paralelos}
S\partial(M,Y)=\{ x \in M\ | \ Y(x) \in T_xM \}.
\end{equation}
\end{defi}

\begin{rema}
\em
The shadow boundary is a natural subset of $M$, it is the
locus where $Y$ is tangent to $M$.\\
The subset $S\partial(M,Y) \subset M$ is closed, so if
$M$ is compact it is also compact. $S\partial(M,Y)$ is not
always a submanifold of $M$. $S\partial(M,Y)$ may be empty
(when $Y$ is nowhere tangent to $M$),
or equal to $M$ (when $Y$ is anywhere tangent to $M$).
\end{rema}

\begin{ejem}
\em
Let $N=S^1 \times S^1 \times S^1$ be the three-torus.
Consider the totally geodesic submanifold $M=S^1
\times S^1 \times \{x\} \subset N$, where $x \in S^1$. Let
$Y: M \longrightarrow TM$ be a parallel vector field on $M$.
Then $Y$ is parallel along $M$. So, $S\partial(M,Y) = M$.\\
Now take $Y' : M \longrightarrow TM^\perp$ orthogonal to $M$.
Again $Y'$ is parallel along $M$ and $S\partial(M,Y') =\emptyset$.\\
Finally, when $N= \mathbb{R}^n$ and $M$ is a compact submanifold,
the shadow boundary with respect to any constant global vector
field $Y$ on $N$ is non-empty.
\end{ejem}

\begin{teor}
\label{sb-productoRiemanniano}
Let $M$ be the Riemannian product $M_1 \times M_2$, of two
submanifolds $M_1 \subset N_1$ and
$M_2 \subset N_2$. Let $Y=(Y_1,Y_2)$ where $Y_j \in
\mathfrak{X}_0(N_j,M_j)$. Then $$S\partial(M,Y)=
S\partial(M_1,Y_1) \times S\partial(M_2,Y_2).$$
\end{teor}
{\bf Proof.}
Let $x=(x_1,x_2) \in S\partial(M,v)$, it means,
that $(Y_1,Y_2) \in T_xM \simeq T_{x_1}M \oplus T_{x_2}M$.
Which is equivalent to $Y_1 \in T_{x_1}M, \ Y_2
\in T_{x_2}M$. Therefore, $x_1 \in S\partial(M_1,Y_1)$ and
$ x_2 \in S\partial(M_1,Y_2)$. This concludes the proof. $\Box$

\begin{ejem}
\label{shadow_of_products}
\em
Let $M=S^1 \times S^1 \subset S^2 \times S^2$.
By Theorem \ref{sb-productoRiemanniano}, we have the following
cases: Let $Y_j \in \mathfrak{X}_0(S^2, S^1)$ (j=1, 2) be
parallel vector fields along $S^1$.\\
Let
      $Y=(Y_1,0)$, with $Y_1 \neq 0$ and $0 \in
      \mathfrak{X}_0(S^2,S^1)$. Let us observe that
      $S\partial(S^1, Y_1)= S^1 \mbox{ or } \emptyset \subset S^2$
      and $S\partial(S^1, 0)= S^1 \subset S^2$.
      Then
      $S \partial(M,Y)= S^1 \times S^1$ or $\emptyset$.\\
      Similarly, if $Y=(0,Y_2)$, $S \partial(M,Y)=S^1
      \times S^1 $ or $\emptyset$.
      Finally $Y=(Y_1,Y_2)$, with $Y_1 \neq 0, \ Y_2 \neq 0$.
      Then $S \partial(M,Y)= S^1 \times S^1$ or $\emptyset$.
\end{ejem}

The second fundamental form of $M \subset N$ at $x \in M$ is a symmetric bilinear tensor, which we denote by
$II_x: T_xM \times T_xM \longrightarrow T_xM^\perp$.
So, $II_x$ is a bilinear application for every $x \in
M$ (see \cite{Ko-No} page 12, for details).

Let $Y$ be a parallel vector field along $M$.
Let $x \in M$ such that $Y(x) \in T_xM$,
then we can consider the following linear application:
$$II(Y(x), \cdot): T_xM \longrightarrow T_xM^\perp.$$
If this transformation is surjective
we will say that $II(Y(x), \cdot)$ is surjective.
In particular, if cod$M=1$, the latter condition is
equivalent to $II(Y(x), \cdot) \neq 0$.

\begin{teor}
\label{Frontera-de-sombra-suave} Let $M$ be a submanifold of
dimension $n$ and codimension $k$ in $N$ ($n \geq k$). Let $Y$ be
a parallel vector field along $M$. If $II(Y(y),\cdot )$ is
surjective for every $y \in S\partial(M,Y)$, then $S\partial(M,Y)$
is a submanifold of dimension $n-k$ in $M$.
\end{teor}
\noindent {\bf Proof.}\\
Let $\nabla$ be the covariant derivative of $N$.
Let $p \in S\partial(M,Y)$, and let $U \subset M$ be a open
neighborhood of $p$. Affirmation: $S\partial(M,Y) \cap U$
is a submanifold of $M$.\\

Let $\xi_j: U \longrightarrow {TU}^\perp$,
$j=1, \ldots, k$ be a basis of orthonormal vector fields
($U$ is such that there exist these vector fields).
Let us consider the next function
$F: U \longrightarrow \mathbb{R}^k$, given by
$$F(x)=(g(Y(x), \xi_1(x)), \ldots, g(Y(x), \xi_k(x))).$$
Is clear that $F^{-1}(0)=S\partial(M,Y) \cap U$. We are going to
prove that $0 \in \mathbb{R}^k$ is a regular value of $F$. We need
verify that for every $x \in S\partial(M,Y) \cap U$, $F_{*x}:T_xM
\longrightarrow \mathbb{R}^k$ is surjective. Let $(y_1, \ldots,
y_n)$ be local coordinates in $U$. Let us calculate the next
derivatives in these coordinates, $\frac{\partial F}{\partial
y_l}= (\frac{\partial}{\partial y_l}g(Y(x), \xi_1(x)), \ldots,
\frac{\partial}{\partial y_l} g(Y(x), \xi_k(x)))$, for every $1
\leq l \leq n$. Since, $Y$ is parallel,
$$\frac{\partial}{\partial y_l}g(Y(x), \xi_j(x))=
g(\nabla_{\partial y_l} Y, \xi_j)+
g(Y, \nabla_{\partial y_l} \xi_j)=
g(Y, \nabla_{\partial y_l} \xi_j).$$
Let us apply Weingarten's formula, which says that\\
$\nabla_{\partial y_l} \xi_j=
-A_{\xi_j}(\partial y_l)+ \nabla^{\perp}_{\partial y_l}\xi_j$.
In conclusion,
$$\frac{\partial}{\partial y_l}g(Y(x), \xi_j(x))=
g( Y , -A_{\xi_j}(\partial y_l) )=
-g(II(Y,\partial y_l),\xi_j(x)),$$ for every
$x \in S\partial(M,Y)$, $1\leq j \leq k$, and
$1 \leq l \leq n$.\\
Now we are ready to see that the next matrix
$$(F_{*x})_{jl}=-(g(II(Y,\partial y_l),\xi_j(x)))$$
has rank $k$.
Let us assume that the row vectors are
linearly dependent,
i.e. we have the following condition\\
$\sum_{j=1}^{k}a_j g(II(Y,\partial y_l),\xi_j(x))=0$, for every
$1 \leq l \leq n$, and where $a_j \in \mathbb{R}$ are constants.
We can rewrite this expression as
$$ g(II(Y,\partial y_l), \sum_{j=1}^{k}a_j\xi_j(x))=0,$$
for every $1 \leq l \leq n$.
Since $II(Y,\cdot) $ is surjective, $\sum_{j=1}^{k}a_j\xi_j(x)=0$,
therefore $a_j=0$. Which proves that $0$ regular value
of $F$. Then we can conclude that, $F^{-1}(0)\cap U$
is a submanifold $U$ of dimension $n-k$. $\Box$

\begin{rema}
\em
A special case of Theorem \ref{Frontera-de-sombra-suave} is
when $\dim N=2 \dim M$. The conclusion in this situation is
that $S\partial(M,Y)$ is a discrete subset of $M$. So if $M$
were compact, $S\partial(M,Y)$ would be a finite set of
points in $M$.
\end{rema}

\begin{rema}
\em In \cite{GRH2}, we proved that when $M^n \subset
\mathbb{R}^{n+1}$ has nowhere zero Gauss-Kronecker curvature,
then for every $v \in \mathbb{R}^{n+1}$, $S\partial(M,v)$ is
a submanifold of $M$ of codimension one.\\
In the particular case of surfaces in $\mathbb{R}^3$, the
smoothness of some shadow boundaries is studied in \cite{Gho},
\cite{Ho} and \cite{Ho2}.
\end{rema}

\begin{defi}
\em
Let $L \subset N$ be a  Riemannian submanifold. Let
$x \in L$, then $L$ is called a {\em totally geodesic
submanifold of $N$ at the point $x$} if every geodesic $\gamma$ of
$L$ through $x$ satisfies $\nabla_{\dot\gamma}
\dot\gamma_{|x}=0$.
\end{defi}
In her work on Affine Differential Geometry \cite{S},
A. Schwenk used conditions similar to those of the next Theorem.
\begin{teor}
\label{ortogonal-implica-tgs}
Let $M^n \subset N^{n+k}$ ($n \geq 2$)
be a submanifold of codimension $k$ ($k \geq 0$).
Let $L \subset M$ be a submanifold of codimension one,
which is not totally geodesic of $N$ at any point.
Let $Y$ be a parallel vector field along $M$ and
orthogonal to $L$. Then
$L \subset S\partial(M,Y)$ if and only if $L$ is a
totally geodesic submanifold of $M$.
\end{teor}
\noindent {\bf Proof.}\\
\noindent $\Longrightarrow$)
Let $x \in L$, since
$\dim(T_xL^{\perp} \cap T_xM)=1$  and by hypothesis
$Y(x) \in T_xL^{\perp} \cap T_xM$, we obtain that
$<Y(x)>=T_xL^{\perp} \cap T_xM$. Therefore,
we have the following equality for every $x\in L$,
\begin{equation}
\label{descomposition} T_xM=T_xL \oplus (T_xL^{\perp} \cap
T_xM)= T_xL \oplus <Y(x)>.
\end{equation}

Let $\gamma \subset L$ be a geodesic and let $x \in \gamma$
be any point.
Hence $\nabla_{\dot{\gamma}} \dot{\gamma}$
is orthogonal to $L$, i.e.
$\nabla_{\dot{\gamma}} \dot{\gamma} \in {T_xL}^\perp.$
Affirmation: $\gamma$ is a geodesic of $M$.
By equality (\ref{descomposition}), we just have to verify that
$\nabla_{\dot{\gamma}} \dot{\gamma}$ is orthogonal to $Y(x)$:
We know that $g(Y(\gamma(t)), \dot{\gamma})=0$, this implies
that, $$g(Y(\gamma(t)), \nabla_{\dot{\gamma}}\dot{\gamma})+
g(\nabla_{\dot{\gamma}}Y(\gamma(t)), \dot{\gamma}) =
\frac{d}{dt}g(Y(\gamma(t)), \dot{\gamma})=0.$$
Then, $\nabla_{\dot{\gamma}} \dot{\gamma}$ is orthogonal
to $M$, so $\gamma$ is a geodesic of $M$.\\
\noindent $\Longleftarrow$) In this implication we assume that
$k=1$. Let $x \in L$, affirmation: $Y(x) \in T_xM$.
Since $L$ is not a totally geodesic submanifold of $N$ at $x$,
there exists a geodesic $\gamma$ of $L$ through $x$ with
$\nabla_{\dot{\gamma}}\dot{\gamma}_{|x} \neq 0$.
By hypothesis, $\gamma$ is also a geodesic
of $M$. So, $\nabla_{\dot{\gamma}}\dot{\gamma} \in
(T_{\gamma}M)^\perp$. Let us prove that $Y(x)$ is orthogonal to
$\nabla_{\dot{\gamma}}\dot{\gamma}$. For this, let us observe that
$g(\dot{\gamma},Y(x))=0$. Therefore,
$$g(Y(\gamma(t)), \nabla_{\dot{\gamma}}\dot{\gamma})+
g(\nabla_{\dot{\gamma}}Y(\gamma(t)), \dot{\gamma}) =
\frac{d}{dt}g(Y(\gamma(t)), \dot{\gamma})=0.$$
But $\nabla_{\dot{\gamma}}Y(\gamma(t))=0$, because $Y$ is
parallel along $L$. Since $M$ is of codimension one,
$g(Y(\gamma(t)), \nabla_{\dot{\gamma}}\dot{\gamma})=0$
implies that $Y(x) \in T_xM$.
$\Box$

\begin{rema}
\em
In Theorem \ref{ortogonal-implica-tgs}, the condition that $L$
is not totally geodesic in $N$ at any point is important to prove
that $L \subset S\partial(M,Y)$. We can see this with the next  example: $N= \mathbb{R}^n$, $M$ a hyperplane, $L$ a linear subspace
of codimension one in $M$. Finally let $Y=v$ be any constant vector
field orthogonal to $M$. In this example the affirmation
$L \subset S\partial(M,Y)$ is false.
\end{rema}

The next Theorem was the original motivation to consider
Helix submanifolds in this work.
\begin{teor}
\label{tgs-implica-helice}
Let $M^n \subset N^{n+k}$ ($n \geq 2$) be
a submanifold of codimension $k$ ($k \geq 0$).
Let $L \subset M$ be a submanifold
and let $Y \in \mathfrak{X}_0(N,M)$.
Assume that $L \subset S\partial(M,Y)$.
If $L$ is a totally geodesic submanifold of $M$,
then $L$ is a helix submanifold of $N$ with respect to $Y$.
\end{teor}
\noindent {\bf Proof.}\\
If $Y(x)\in T_xL$, for every $x\in L$, $L$ is a helix.
Otherwise, let $p\in L$ such that $Y(p) \notin T_pL$.
So,
\begin{equation}
\label{orthogonal-cap3}
T_pL \oplus <Y(p)> \subset T_pM.
\end{equation}

Affirmation:
The angle between $T_xL$ and $Y(x)$ is constant,
for every $x$ in $L$.
Let $\gamma$ be any geodesic of $L$ from $p$ to $x$,
hence it is also geodesic of $M$.
Now, let us consider the parallel transport
$\tau$ in $M$, along $\gamma$, from  $p$ to $x$. Therefore,
$\tau : T_pM \longrightarrow T_xM$
is an isometry. So, $\tau$ transforms the latter
equation (\ref{orthogonal-cap3}), in
$T_xL \oplus <Y(x)> \subset T_xM$.
Since the parallel transport is an isometry,
the angle between $T_xL$ and $Y(x)$ is
equal to the angle between
$T_pL$ and $Y(p)$. $\Box$

\begin{rema}
\em
Let us observe that any Euclidean compact helix with respect to
$Y=v \in \mathbb{R}^{n+k} -\{0\}$ should be orthogonal to $Y$.
So, if $N=\mathbb{R}^{n+k}$ and $L \subset S\partial(M,Y) \subset M$
is a compact totally geodesic submanifold of $M$, then
by Theorem \ref{tgs-implica-helice}, $L$ is orthogonal to $Y$.
\end{rema}

\end{section}

\begin{section}{Minimal shadow boundaries}

\begin{defi}
\em
Let $L \subset M$ be a Riemannian submanifold.
Let $x \in L$, $n=\mbox{dim}L$ and let $e_1, \ldots ,e_n$
be an orthonormal basis of $T_xL$. The
{\em mean curvature vector field} $H$ of $L \subset M$
at $x$, is
$$H(x)=\frac{1}{n} \sum_{i=1}^n II_x(e_i,e_i).$$
\end{defi}
\begin{rema}
\em
For every $x \in L$, $H(x) \in T_xL^\perp$, because
$II_x(e_i,e_i) \in T_xL^\perp$ where $i=1,\ldots,n$.
\end{rema}
\begin{defi}
\label{minimal}
\em
A submanifold $L \subset M$ is {\em minimal} if
$H(x)=0$ for every $x \in L$.
\end{defi}
The next proposition and its proof is due to C. Bang-yen,
see \cite{Ba-yen}.
\begin{prop}{(Bang's Lemma)}
\label{lemadeBang}
Let $L^n$ be a submanifold of $M^s$, where $M$ is a submanifold
of $N^m$. Then $L$ is minimal in $M$
if and only if the mean curvature vector field of
$L \subset N$ is orthogonal to $M$.
\end{prop}
{\bf Proof.}
Let $X$ and $Y$ be two vector fields on $L$. Let
$\nabla$ and $\nabla'$ be the covariant derivatives
of $N$ and $M$ respectively.
Gauss formula for $M \subset N$ says that
$$\nabla_X Y = \nabla'_X Y+II^N(X,Y),$$ where $II^N$ is the
second fundamental form of $M \subset N$.\\
Let $\nabla''$ be the covariant derivative of $L$ and $II^M$ the second
fundamental form of $L \subset M$. Then we have
$$\nabla'_X Y = \nabla''_X Y+II^M(X,Y).$$
From the two latter formulae we get
$$\nabla_X Y = \nabla''_X Y+II^M(X,Y)+II^N(X,Y).$$
So, the second fundamental form $II$ of
$L \subset N$ is $$II(X,Y)=II^M(X,Y)+II^N(X,Y).$$
By definition, $II^M$ is orthogonal to $L$, tangent to $M$ and
$II^N$ is orthogonal to $M$. Now, let us consider the mean curvature vector fields.
Let $H$ and $H^M$ be the mean curvature vector fields
of $L \subset N$ and $L \subset M$ respectively.
Let $x \in L$ and let $e_1, \ldots , e_n$ be a orthonormal
of $T_xL$. Then
$\frac{1}{n} \sum_{i=1}^n II(e_i,e_i)=
\frac{1}{n} \sum_{i=1}^n II^M(e_i,e_i)+
\frac{1}{n} \sum_{i=1}^n II^N(e_i,e_i)$.
Hence,
$$H(x)=H^M(x)+H(L;M,N),$$
where $H(L;M,N)=\frac{1}{n} \sum_{i=1}^n II^N(e_i,e_i)$.
Let us observe that $H(L;M,N)$ is orthogonal to $M$.
Then $L$ is minimal in $M$ if and only if
$H(x)=H(L;M,N)(x)$.$\Box$\\

\begin{teor}
\label{minimal-contenida-fs}
Let $M^n \subset N^{n+1}$ be a Riemannian submanifold
and let $L^{n-1} \subset M$ be a submanifold such that
$L \subset S\partial(M,Y)$, where $Y$ is parallel along
$M$ and transverse to $L$. Let $H$ be the mean curvature
vector field of $L \subset N$.
Then $L$ is minimal in $M$ if and only if $g(H,Y)=0$.
\end{teor}
{\bf Proof.}
By hypothesis $Y(x) \in T_xM$ for every $x \in L$.
By Lemma \ref{lemadeBang}, if $L$ is minimal in $M$ then $H$
is orthogonal to $M$. So, $H(x)$ is orthogonal to $Y(x)$,
i.e. $g(H,Y)=0$.\\
Now, let us assume that $g(H,Y)=0$.
By definition, $H$ is orthogonal to $L$. To apply Proposition
\ref{lemadeBang}, we need prove that $H$ is orthogonal to $M$.
Since $Y$ is transversal to $L$, $T_xM=T_xL \oplus <Y(x)>$
for every $x \in L$. Now it is clear that $H$ is orthogonal to $M$.
Then $L$ is minimal in $M$.
$\Box$

\begin{coro}
Let $N = \mathbb{R}^{n+1}$, let $M^n \subset N$ be a
submanifold and let $Y$ be a constant vector field on $N$
induced by $v \in \mathbb{R}^{n+1} - \{0\}$.
Let $L \subset S\partial(M,Y)$ be a compact submanifold
of $M$ of codimension one. Let us assume that
\begin{itemize}
\item $L$ is contained in a hyperplane,
which is transverse to $v$.
\end{itemize}
If $L$ is minimal then $L$ is a totally
geodesic submanifold of $M$.
\end{coro}
\noindent {\bf Proof.}
Let $\Pi$ be the hyperplane that contains $L$. Hence the
mean curvature vector $H$, of $L \subset \mathbb{R}^{n+1}$,
is contained in $\Pi$: $H_p \in T_p\Pi$
for every $p \in L$. Since $\Pi$ is transversal to $v$,
then $v$ is not tangent to $L$.
Now we can apply Theorem \ref{minimal-contenida-fs}
($M$ is minimal), to deduce that $<H(p),v>=0$, for every
$p \in L$. In summary, $v$ is orthogonal to $H(p)$, where
$H(p) \in T_p\Pi$.\\
Affirmation: $v$ is orthogonal to $\Pi$. Since $L$ is compact,
we can get a basis of $T_p \Pi$, by translation of vectors
$H(x)$ into the point $p$.
Moreover,
$H(x)$ can not be zero for every $x$ because $L$ is compact
and $N$ do not contains minimal compact submanifolds.\\
We proved that $H(x)$ is orthogonal to $v$, for every $x$
in $L$. This proves that $v$ is orthogonal to $\Pi$. Since
$L \subset \Pi$, $v$ is orthogonal to $L$.\\
Finally, we can apply Theorem \ref{ortogonal-implica-tgs},
which says that if $L \subset S\partial(M,Y)$ and
$v$ is orthogonal to $L$, then $L$ is a totally geodesic
submanifold of $M$.
$\Box$

\begin{ejem}
\label{minimal-contenida-fs}
\em
In this example we are going to construct a submanifold $M$ of
$N=\mathbb{R}^{n+2}$, which contains a minimal submanifold
in some shadow boundary.\\
Let $L^{n} \subset \mathbb{R}^{n+2}$ be a submanifold
and let $Y=v \in \mathbb{R}^{n+2}-\{0\}$ such that:
\begin{itemize}
\item $v$ is transverse to $L$: $v \notin T_xL$
for every $x \in L$,
\item $<H,v>=0$,
where $H$ is the mean curvature vector field
of $L \subset N$.
\item
$L_{\epsilon,v}=\{ y= x+\lambda v \in \mathbb{R}^{n+2} \ | \ x
\in L, \ |\lambda| < \epsilon \}$ is a submanifold, where
$\epsilon$ denotes a positive smooth function of $L$.
\end{itemize}
Then $L$ is a minimal submanifold of $M=L_{\epsilon,v}$.
If $L$ is compact, $\epsilon$ can be a constant function.
Proof: this is consequence of Theorem \ref{minimal-contenida-fs}.
We should verify the hypothesis of such theorem.
The submanifold $M=L_{\epsilon,v}$ is a
``cylindric" neighborhood of $L$ in direction $v$.
For this reason,
$S\partial(M, v)=L_{\epsilon,v}=M$,
and then $L \subset S \partial(M, v)$.
\end{ejem}

\begin{ejem}
\em

Construction in Example \ref{minimal-contenida-fs} is based in a
submanifold of codimension two, which is transverse to $v$ and
whose mean curvature vector field is orthogonal to $v$. There are
submanifolds of codimension two which does not admits a transverse
direction like
$S^1 \times S^1 \subset \mathbb{R}^4$.\\
There are examples which satisfy conditions of Example
\ref{minimal-contenida-fs}:
\begin{enumerate}
\item
Let $N_1 \subset N= \mathbb{R}^3$ be a classic helix with respect
to direction $v$, in fact it can be a planar curve orthogonal
to $v$.
Let us assume that its curvature is nowhere zero.
Since $N_1 \subset N$ is a helix, its mean curvature vector
field, $H_1$, is orthogonal to $v$.
In fact $||H_1|| \neq 0$ is its curvature.\\
Let us consider the surface $M=N_1 \times \mathbb{R} \subset
\mathbb{R}^3 \times \mathbb{R}=\mathbb{R}^4$.\\
Affirmation: $<H_M,v>=0$ and $v \in \mathbb{R}^3 \times \{0\}$ is
transverse to $M$, where $H_M$ is the mean curvature vector
of $M \subset \mathbb{R}^4$:\\
$H_M(x,y)=\frac{1}{2}H_1(x)$. Since $v$ is transverse to $\{0\}
\times \mathbb{R}$ and to $N_1$, it is transverse to $M$.
\item
It is not necessary to take a classic helix in the latter example.
Here is another construction. Let $N_1 \subset \mathbb{R}^2$ be
any embedded curve whose curvature is nowhere zero. Let $N_2
\subset \mathbb{R}^3$ be any minimal surface of $\mathbb{R}^3$
which admits a transverse direction $w \in \mathbb{R}^3$. Let
$M^3=N_1 \times N_2 \subset \mathbb{R}^2 \times \mathbb{R}^3
=\mathbb{R}^5$. Affirmation: $w \in \{0\} \times \mathbb{R}^3$ is
transverse to $M$ and
$<H_M,w>=0$.\\
As before, $H_M(x,y)=\frac{1}{3}H_1(x) + \frac{2}{3}
H_2(y)=\frac{1}{3}H_1(x)$. Now let us observe that $H_1(x) \in
\mathbb{R}^2 \times \{0\}$, $w$ is orthogonal to $N_1$ and
transverse to $N_2$.
\end{enumerate}
\end{ejem}

\begin{question}
\em
Does exist a compact strictly convex hypersurface in
$\mathbb{R}^{n+1}$ with some minimal and non-totally
geodesic shadow boundary?
\end{question}

\end{section}

{\bf Acknowledgements}
This paper is based on part of the author's Ph.D. thesis.
\cite{GRH}. He wishes to express his gratitude to
Luis Hern\'andez Lamoneda, his thesis advisor,
for his support and helpful suggestions.\\

\noindent {\bf Gabriel Ruiz-Hern\'andez }\\
{\footnotesize CIMAT, Jalisco s.n., Mineral de
Valenciana Guanajuato 36240, M\'exico;\\
\noindent e-mail: gruiz@cimat.mx}\\

\end{document}